\documentclass[12pt, amssymb]{article}

\usepackage{latexsym}   
\usepackage{amssymb}    
\usepackage{amsmath}    
\usepackage{amsthm}     

\newtheorem{theorem}{Theorem}[section]
\newtheorem{lemma}[theorem]{Lemma}
\newtheorem{definition}[theorem]{Definition}
\newtheorem{corollary}[theorem]{Corollary}
\newtheorem{proposition}[theorem]{Proposition}

\newtheorem{remark}[theorem]{Remark}

\def\X{{\cal X}}
\def\Y{{\cal Y}}
\def\O{{\cal O}}
\def\D{{\cal D}}

\def\E{{\cal E}}
\def\F{{\mathbb F}}
\def\Z{{\mathbb Z}}
\def\G{{\mathbb G}}
\def\Q{{\mathbb Q}}
\def\C{{\mathbb C}}
\def\c{{\mathcal C}}
\def\R{{\mathbb R}}

\begin{document} 

\title{$\epsilon$-Constants and Orthogonal Representations} 

\author{Darren Glass}\date{} 

\maketitle

\begin{abstract}
In this paper we suppose $G$ is a finite group acting tamely on a regular
projective curve $\X$ over $\Z$ and $V$ is an orthogonal representation of $G$
of dimension $0$ and trivial determinant.  Our main result determines the sign
of the $\epsilon$-constant $\epsilon(\X/G,V)$ in terms of data associated to
the archimedean place and to the crossing points of irreducible components of
finite fibers of $\X$, subject to certain standard hypotheses about these
fibers.
\end{abstract}

\section{Introduction}
\label{section:intro}
\setcounter{equation}{0}

This section will state the main questions and results of the paper, specify
notation, and give some background. Let $\X$ be an arithmetic scheme of
dimension $d+1$ which is flat, regular, and projective over $\Z$.  We suppose
that $f:\X \rightarrow Spec(\Z)$ is the structure morphism and that its fibres
are all of dimension $d$.  Let $G$ be a finite group which acts tamely on $\X$
in the sense that for each closed point $x \in X$, the order of the inertia
group of $x$ is relatively prime to the residue characteristic of $x$. Define
$\Y$ to be the quotient scheme $\X / G$. We assume that $\Y$ is regular, and
that for all finite places $v$  the fiber $\Y_v = (\X_v)/G = \Y
\otimes_\Z(\Z/p(v))$ has normal crossings and smooth irreducible components
with multiplicities relatively prime to the residue characteristic of $v$.
Finally, let $V$ be a representation of $G$ over $\overline{\Q}$.

Associated to this data, there are well-known $\zeta$-functions and
$L$-functions, both functions of a complex variable $s$. The $L$-function is
conjectured to have a functional equation of the form:  $$L(s,\Y,V) =
\epsilon(\Y,V)A(\Y,V)^{-s}L(d+1-s,\Y,V^*)$$

\noindent in which $A(\Y,V)$ is a positive integer called the conductor, $V^*$
is the dual representation of $V$ and the $\epsilon$-constant $\epsilon(\Y,V)$
is a nonzero algebraic number. In recent years, many authors have studied the
problem of determining these $\epsilon$ constants, which may be defined
unconditionally after we choose an auxiliary prime $\ell$. 

This paper concerns the case where $V$ is an orthogonal representation, meaning
that there is a non-degenerate symmetric $G$-invariant bilinear form $V \times
V \rightarrow \overline{\Q} \subseteq \C$ (where we fix an embedding of
$\overline{\Q}$ into $\C$). In order to get the strongest results, we will
furthermore make the technical hypotheses that $V$ is a virtual representation
of trivial determinant and dimension zero.  In other words, $V$ will be a
linear combination of orthogonal representations such that the weighted sum of
their dimensions is zero and the product of their determinants is trivial.

We can now state in general terms the main result of this paper.

\begin{theorem}
\label{thm:main}

If $d=1$ and $V$ is an orthogonal virtual representation of degree zero
and trivial determinant then the sign of the constant $\epsilon(\Y,V)
\in \R^*$ can be determined from the $\epsilon$-constant
$\epsilon_\infty(\Y,V)$ and from the restriction of the $G$-cover $\X
\to \Y$ over the finite set of closed points $z$ of $\Y$ where two
distinct irreducible components of a fiber of $\Y$ over $Spec(\Z)$
intersect.

\end{theorem}

The constant $\epsilon_\infty(\Y,V)$ which comes up in this formulae is the
archimedean $\epsilon$-constant defined by Deligne in \S 8 of \cite{D1} using
the action of the group $G$ and of complex conjugation on the Hodge cohomology
groups $H^{p,q}(\X,\C)$. Section Two of this paper recalls this and other
definitions of $\epsilon$-constants, as well as work done by Deligne, Frohilch,
Queyrut, Chinburg, Erez, Pappas, and Taylor in computing $\epsilon$-constants
associated to situations similar to those in Theorem \ref{thm:main}.  In
Section \ref{section:main}, we make a more precise statement of the main
theorem and prove it.  The proof uses formulae of Saito, Classfield theory, and
several of the results discussed in Section \ref{section:bkgrnd}.

The main results of this paper are from the author's dissertation, and he would like to express his gratitude to his advisor, Ted Chinburg.  

\section{Background}
\label{section:bkgrnd}
\setcounter{equation}{0}

In this section, we look at some of the work that others have done in
order to compute $\epsilon$-constants in various situations. 

Fr\"ohlich and Queyrut look at computing $\epsilon$-constants in the
case where $\X$ and $\Y$ are of relative dimension $0$ over $\Z$ and
$V$ is an orthogonal representation.  In {\cite{FQ}} they are able to
prove the following result:

\begin{theorem}
\label{thm:fq}

If $d = 0$ and $V$ is an orthogonal representation of $G$, then
$\epsilon(\Y,V)$ is positive.

\end{theorem}

We now recall some elements of Deligne's theory of local constants, which is essential for our work.

\begin{definition} Let $\X$ and $\Y = \X/G$ be as above. Let $V$ be any
virtual complex representation of $G$.

\begin{enumerate}

\item[a.] Let $\epsilon_{v,0}(\Y,V)$ be the Deligne local constant
defined in \cite{CEPT1}. (see also \cite{D1}). In particular, the
definition of $\epsilon_{v,0}(\Y,V)$ requires that one chooses an
auxiliary prime $\ell \neq v$, a nontrivial continuous complex
character of $\Q_v$ which we denote by $\psi_v$ and a Haar measure
$dx_v$ on $\Q_v$.  In the case where $V$ has trivial determinant and is
of dimension $0$, then $\epsilon_{v,0} (\Y,V)$ is independent of these
choices (see Proposition 2.4.1 of \cite{CEPT1}).  This term is
well-defined for $v = \infty$ as well as for finite places $v$.

\item[b.] Let $X$ be a variety of dimension $d$ which is defined over a
finite field of characteristic $p$. Let $\ell$ be a prime different
from $p$ and let $j_\ell : \overline{\Q_\ell} \rightarrow \C$ be an
embedding.  Finally, define $V_\ell$ to be a virtual representation of
$G$ over $\overline{\Q_\ell}$ such that $j_\ell(\chi_{V_\ell}) =
\chi_V$.  Define $\epsilon(X,V) = j_\ell(det(-F | (H^*_{et}
(\overline{\F_p}\times_{\F_p}X,\Q_{\ell}) \otimes V_\ell^*)^G))$, where
$F$ is the geometric Frobenius automorphism.  This number is
independent of all choices.

\item[c.] For finite places $v$ of $\Q$, we let $\epsilon_{v}(\Y,V) =
\epsilon_{v,0}(\Y,V) \epsilon(\Y_v,V)$, where $\epsilon(\Y_v,V)$ is
defined as an $\epsilon$-constant over a finite field.  Furthermore,in
the case where $v = \infty$, we let $\epsilon(\Y_v,V) = 1$ so that in
particular $\epsilon_\infty(\Y,V) = \epsilon_{\infty,0}(\Y,V)$.

\item[d.] The global $\epsilon$-constant associated to $V$ is defined by
$\epsilon(\Y,V) = \prod_v \epsilon_v(\Y,V)$ where the product is over all
places $v$ of $\Q$.

\end{enumerate}
\end{definition}

The $\epsilon$-constants associated to varieties defined over finite
fields are studied by Chinburg, Erez, Pappas, and Taylor in
{\cite{CEPT3}}.  Other papers by these authors, such as {\cite{CEPT1}} and
{\cite{CEPT2}} prove results on computing $\epsilon$-constants
associated to arithmetic schemes in the case where $V$ is a symplectic
representation.  Recall that a symplectic representation is a representation $V$ which is equipped with a non-degenerate alternating $G$-invariant bilinear form. 

For many applications of $\epsilon$-constants it is not the
actual $\epsilon$-constant we are interested in computing but merely
the sign of this constant.  We denote the sign of $\epsilon(\Y,V)$ by $W(\Y,V)$ and
call this the root number of $V$.

\section{Main Results}
\label{section:main}
\setcounter{equation}{0}

\subsection{Reduction To Fibral Computations}

Let $\X$, $G$, $\Y = \X/G$ be as in \S \ref{section:intro}.  Let $S$ be
the set of all finite places $v$ of $\Q$ where either the fiber $\Y_v =
\Y \otimes_\Z(\Z/p(v))$ is not smooth or the map $\pi:\X \rightarrow
\Y$ is ramified.  Let $D'$ be a horizontal divisor on $\Y$ such that
$D'+ \Y_T = K_{\Y} + \Y_S^{red}$, where $K_{\Y}$ is a canonical divisor
on $\Y$, $\Y_S^{red}$ is the sum of the reductions of the fibers of
$\Y$ at the places in $S$, $T$ is a finite set of finite places of $\Q$
which is disjoint from $S$, and $\Y_T$ is the sum of the (necessarily
reduced) fibers of $\Y$ over the places in $T$.  Thus $O_\Y(D' + \Y_T)$
is isomorphic to the twist $\omega_{\Y/\Z}(\Y_S^{red})$ of the relative
dualizing sheaf $\omega_{\Y/\Z}$ by $\O_\Y(\Y_S^{red})$.  We  further
wish to choose $D'$ so that it intersects the non-smooth fibers $\Y_v$
of $\Y$ transversally at smooth points on the reduction of $\Y_v$.  We
can choose such a $D'$ after a suitable base change due to the moving
lemma proven as Proposition $9.1.3$ in \cite{CEPT1}. The choice of this
canonical divisor is not unique, but our calculation will show that the
results are independent of the choice of $D'$.   

As stated above, we can only choose a horizontal divisor $D'$ with the desired
properties after a suitable base change.  Thus, we need to consider how base
changes will affect the $\epsilon$-constants.  To be precise about how we make
the base change, we will choose an odd prime $\ell$ which is not in the set of
bad primes $S$, and we denote by $N_\infty$ the cyclotomic $\Z_\ell$ extension
of $\Q$.  Because we have chosen $\ell \notin S$, this base extension is
\'etale over $S$, and the pullback of a canonical divisor remains canonical up
to a multiple of the fiber of $\Y$ over $\ell$.   Proposition $9.1.3$ of
\cite{CEPT1} shows that a horizontal divisor $D'$ with the required properties
exists after a base extension to a the ring of integers of a finite extension
of $\Q$ inside $N_\infty$.  This base extension, which we now fix, is of degree
a power of $\ell$.  Since $\ell$ is not in the set $S$, the Hasse-Davenport
Theorem together  with Lemma $9.4.1$ of [CEPT1] shows that the epsilon
constants we will consider for the base change are the $\ell^a$-th power of the
corresponding constants before the base change.  Because we are primarily
interested in the sign of the $\epsilon$-constant, we are free to make a base
change of the above kind.   If we were interested in preserving more
information about $\epsilon$ we will place a stricter congruence condition on
the prime $\ell$. 

\begin{lemma}
For the infinite place, $\epsilon_{\infty,0}(D',V) = 1$
\end{lemma}

This lemma is an immediate corollary to Proposition 5.4.2 of
\cite{CEPT1}.  In particular, this proposition says that if $d$ is odd then the
archimedean epsilon constant associated to the canonical divisor and to
any representation $V$ of trivial determinant and dimension zero is equal
to one.  $D'$ differs from the canonical divisor only by vertical
fibers, and thus the result applies.

\begin{lemma}
\label{lemma:epvo}
With $\Y$, $D'$, and $V$ chosen as above, $\epsilon_{v,0}(\Y,V) =
\epsilon_{v,0} (D',V)$ for all finite places $v$ of $\Q$.

\end{lemma}

\noindent {\bf Proof:}\enspace For all places $v \in S$, this follows directly
from \cite{CEPT1}.  However, one can generalize their results in order to prove
the lemma.  In order to do this, we let $\c_v$ be the set of irreducible
components of $\Y_v^{red}$.  For each $C_i \in \c_v$ let $\kappa_i$ be the
Gauss sum associated to the restriction of the representation $V$ to the
inertia group of the generic point of $C_i$ as defined in \cite{CEPT1}.  They
define $c_i$ to be the $\ell$-adic Euler characteristic with compact support of
the open subscheme of $C_i$ consisting of points which are nonsingular in
$\Y_v^{red}$.  The formulae developed by Saito as Theorems $1$ and $2$ of 
\cite{tS} imply that $\epsilon_{v,0}(\X,V) = \prod_{i \in \c_v}
\kappa_i(V)^{c_i}$.

For each $C_i$ we compute that
$deg_{C_i}(\O_\Y(K_\Y+\Y_S^{red})) = -c_i f_i$, where $f_i$ is the
index of the constant field extension $[F_i:\F_p]$. Changing views, we
let $\delta'$ be a point where $\Y_v^{red}$ intersects the horizontal
divisor $D'$.  We define Gauss sums $\kappa_{\delta'}$ in a similar
way to the above defined $\kappa_i$, such that, in particular,
$\kappa_{\delta'} = \kappa_i^{[k(\delta):F_i]}$.  Furthermore, the
local epsilon constant $\epsilon_{v,0}(D',V)$ is given by
$\prod_{\delta' \in D' \cap \Y_v^{red}}\kappa_{\delta'}$ (see \cite{tS}
p. 416). The proof of the lemma in this case now reduces to counting
intersection numbers and verifying that $\kappa_i$ occurs as a factor
the same number of times in both  $\epsilon_{v,0}(\Y,V)$ and
$\epsilon_{v,0} (D',V)$.

For the finite places $v$ which are not in $S$ the argument is
similar.  It is only the intersection multiplicities of $D'$ with
certain vertical divisors that matters, and these numbers do not change
in the event that we add new vertical fibers into the divisors. For this
reason, the appearance of $\Y_T$ in the equality $D'+ \Y_T = K_{\Y} +
\Y_S^{red}$ makes no difference in the argument.   $\blacksquare$

With these lemmas in hand, we can make the following series of calculations:
\begin{eqnarray}
\label{eq:doit}
\epsilon(\Y,V) &=& \prod_v\epsilon_v(\Y,V) \nonumber \\
&=&\epsilon_\infty(\Y,V) \prod_{v finite} \epsilon_{v,0}(\Y,V)\epsilon(\Y_v,V) \nonumber\\
&=&\epsilon_\infty(\Y,V)\prod_{v finite} \epsilon_{v,0}(D',V)\epsilon(\Y_v,V) \nonumber\\
&=&\epsilon_\infty(\Y,V)\epsilon_{\infty,0}(D',V)\prod_{v finite}\epsilon_{v,0}(D',V)\epsilon(D'_v,V)\epsilon(D'_v,V)^{-1}\epsilon(\Y_v,V)\nonumber\\
&=&\epsilon(D',V)\epsilon_{\infty,0}(\Y,V) \prod_{v finite} \epsilon(D'_v,V)^{-1}\epsilon(\Y_v,V)
\end{eqnarray}

In these calculations, $D'_v = D' \otimes_\Z \Z/p(v)$ is the finite
collection of closed points of $D'$ lying above the finite place $v$ of $\Q$.

\begin{lemma}
\label{lemma:fq}
$\epsilon(D',V)$ is positive.
\end{lemma}

\noindent {\bf Proof:}\enspace $D'$ is a one-dimensional object, and
the restriction of $V$ to $D'$ will still be an orthogonal representation.
By applying the theorem of Fr\"ohlich-Queyrut to the normalization of
$D'$ (which we denote by $(D')^\#$), we get that $\epsilon((D')^\#,V)$
is positive. Now, because the definition of local constants involves
only the Galois action on general fibers, $\epsilon_{v,0}((D')^\#,V) =
\epsilon_{v,0}(D',V)$.  Thus, we are only concerned with the difference
between the terms $\epsilon((D')_v^\#,V)$ and $\epsilon(D'_v,V)$, all
of which come about from the singular points $z$ of $D'$.  The action
of $G$ is \'etale at these points, and thus we can compute the local
constants at these points as $\epsilon(y,V)=det(-F | (H^0(y,\Q_\ell)
\otimes V)^{G_x}=det(V)(\pi_{\Y_v^{red},y})$, which is equal to
one due to our hypotheses that $V$ has trivial determinant. 
$\blacksquare$

Thus, we have reduced the calculation of the sign of $\epsilon(\Y,V)$,
which is an inherently two-dimensional calculation, to a collection of
fibral computations $\epsilon(D'_v,V)^{-1}\epsilon(\Y_v,V)$
for each finite place $v$, and a calculation for the archimedean component
$\epsilon_{\infty,0}(\Y,V)$.

\subsection{The One-Component Case}

\begin{theorem}
\label{thm:irred}

Let $\X, \Y, D'$ be as above and let $V$ be an orthogonal representation of
trivial determinant and dimension.  Furthermore, assume $v$ is a finite place
of $\Q$ such that $\Y_v^{red}$ is irreducible. Then
$\epsilon(D'_v,V)^{-1}\epsilon(\Y_v,V) = 1$.

\end{theorem}

{\bf Proof:} Assume that $\Y_v^{red}$ consists of a single component.  Then
$\Y_v^{red}$ is smooth by hypothesis. Let $c$ be an irreducible component of
$\X_v$ with generic point $\mu_c$. Let $G_{\mu_c}$ be the Galois group acting
on the generic point of $c$, and $I_{\mu_c}$ be the inertia group at the
generic point of $c$. Then we have that $I_{\mu_c} \subseteq G_{\mu_c}
\subseteq G$. We denote $I_{\mu_c}$ by $I$.   The tameness hypotheses implies
that the order of $I$ is relatively prime to $v$, and that $I$ is a cyclic
group.  The specific structure of $I$ is discussed in detail in the Appendix to
\cite{CEPT1}.

We begin by computing $\epsilon(\Y_v,V) = \prod_i
det(-F|(H^i(\overline{\Z/v\Z}$ $\otimes_{\Z/v\Z}$ $\Y_v,\Q_\ell)\otimes
V)^G)^{(-1)^{i+1}}$, where $F$ is the Frobenius element as described
above. We know by our hypotheses that the cover $\X_v^{red}\rightarrow
\Y_v^{red}$ is a tame $G_{\mu_c}/I$-cover of smooth curves over
$\Z/p\Z$. Furthermore, the action of $G/I$ on $\X_v^{red}$ is \'etale
because $I = I_{\mu_c} = I_{\X,x}$, the inertia group of the point $x$, for
all points $x \in \X_v^{red}$.   This implies that  $\epsilon(\Y_v,V) =
\epsilon(\Y_v^{red}, V^{I})$.  $I$ acts trivially on the cohomology
group $H^*(\X_v^{red},\Q_\ell)$, so Saito's formulae in \cite{tS} imply
that $\epsilon(\Y_v,V)$ can be calculated as
$det(V^I)(K_{\Y_v^{red}})$, where $K_{\Y_v^{red}}$ is the canonical
divisor on $\Y_v^{red}$.  The terms $K_{\Y_v^{red}}$ are well defined,
as we have assumed that for all finite $v$ the irreducible components of
$\Y_v^{red}$ are themselves smooth.

Next we look at the term $\epsilon(D_v',V)$. Let $D$ be the preimage of
$D'$ in $\X$, and let $\ell$ be a prime different from $v$. Let
$I_{D,x}$ be the cyclic inertia group of a point $x$ lying above the
points in $\Y_v \cap D'$ (note that this is independent of which point
$x$ we choose). Because $D_v'$ is zero dimensional, we know
that $\epsilon(D_v',V) = \prod_{y \in D_v'} \epsilon(y,V)$ where $$\epsilon(y,V) =
det(-F|(H^0(\pi^{-1}(y)^{red},\Q_\ell)\otimes V)^G)$$ if we view $\pi$ as
the cover $D_v \rightarrow D_v'$.  We know that $\pi^{-1}(y)^{red} = (y
\times_{D'_v}D_v)^{red} = x \times_{G_x} G$. In particular, this implies that

$$H^0(\pi^{-1}(y)^{red},\Q_\ell)\otimes V = (Ind^G_{G_x} H^0(x,
\Q_\ell)) \otimes V$$ Because $I_{\mu_c} = I_{\X,x}$ for all points
$x$, we can see that  $$Ind^G_{G_x} H^0(x,\Q_\ell) =
Infl^G_{G/I_{\X,x}} Ind^{G/I_{\X,x}}_{G_x/I_{\X,x}} H^0(x,\Q_\ell)$$ 
Recall that $I = I_{\X,x}$ acts trivially on $H^0(x, \Q_\ell)$. This
allows us to compute that 

\begin{eqnarray*}
\epsilon(y,V)&=&det(-F|(H^0(\pi^{-1}(y)^{red},\Q_\ell)\otimes V)^G) \\
&=&det(-F|(Infl^G_{G/I}
Ind^{G/I}_{G_x/I} H^0(x,\Q_\ell)\otimes V)^G)\\
&=&det(-F| (Ind^{G/I}_{G_x/I} H^0(x, \Q_\ell) \otimes V^{I})^{G/I})\\
&=&\epsilon(y,V^I)
\end{eqnarray*}

\noindent where $\epsilon(y,V^I)$ is the local constant associated to the $G/I$
cover $\X_v^{red} \rightarrow \Y_v^{red}$.

This last term is in turn equal to $det(V^I)(\pi_{\Y_v^{red},y})$,
where $\pi_{\Y_v^{red},y}$ is the local uniformizer from classfield
theory since $\X_v^{red} \rightarrow \Y_v^{red}$ is an unramified $G/I$
cover. Finally, we can put these terms together to get that $\epsilon(D'_v,V) =
det(V^I)(D' \cap \Y_v^{red})$, where $D' \cap \Y_v^{red}$ is viewed as
a divisor on $\Y_v^{red}$.

\begin{lemma}
\label{lemma:canonical}

Under the above hypotheses, $D' \cap \Y_v^{red}$ is a canonical divisor
on $\Y_v^{red}$.

\end{lemma}

Given this lemma, we will have shown that $\epsilon(D'_v,V) = det(V^I)(K) = \epsilon(\Y_v,V)$, so in particular $\epsilon(D'_v,V)^{-1}\epsilon(\Y_v,V) = 1$, and Theorem \ref{thm:irred} will be proven.

In order to prove Lemma \ref{lemma:canonical}, recall that we chose
$D'$ so that $\O_\Y(D'+ \Y_T) = \omega_{\Y/\Z}(\Y_S^{red})$.   We note
that if we look at the two exact sequences: $$ 0 \rightarrow
\O_\Y(-\Y_v^{red}) \rightarrow \O_\Y \rightarrow \O_{\Y_v^{red}}
\rightarrow 0$$ $$ 0 \rightarrow \O_\Y(D' - \Y_v^{red}) \rightarrow
\O_\Y(D') \rightarrow \O_\Y(D')|_{\Y_v^{red}} \rightarrow 0$$ we get
that for all primes $v$, $\O_\Y(D')|_{\Y_v^{red}}$ is the same as
$\O_{\Y_v^{red}}(D' \cap \Y_v^{red})$.  Furthermore, for those primes
$v$ which are in $S$ (and in particular are not in $T$), we further get
that $\O_\Y(D')|_{\Y_v^{red}} = \O_\Y(D' +\Y_T)|_{\Y_v^{red}}$.  We now
are able to make the following computation for all $v \in S$ such that
$\Y_v^{red}$ is irreducible:
\begin{eqnarray*}
\O_{\Y_v^{red}} (D' \cap \Y_v^{red})&=& \O_\Y(D') |_{\Y_v^{red}}\\
&=& \O_\Y(D'+\Y_T) |_{\Y_v^{red}}\\
&=&\omega_{\Y/{\bf Z}}(\Y_S^{red})  |_{\Y_v^{red}} \\
&=&\omega_{\Y/{\bf Z}}(\Y_v^{red})  |_{\Y_v^{red}} \\
&=&\omega_{\Y_v^{red}} \\
\end{eqnarray*}
In other words, for such $v$, $D' \cap \Y_v^{red}$ is a canonical
divisor on $\Y_v^{red}$ under these assumptions. 

It remains to show that Lemma \ref{lemma:canonical} holds for primes
$w$ outside of the set $S$.   We know that for such $w$, the fibres
$\Y_w$ are reduced and smooth and that the local equations have a nice
form.  This implies in particular that $\Y_w^{red} = \Y_w$ is a
principal divisor and thus that $\O_\Y(\Y_w^{red})$ is isomorphic to
$\O_\Y$.

Recall that by definition we have that $D' + \Y_T = K_\Y+\Y_S^{red}$. 
This tells us that
$$D' + \Y_T - \Y_S^{red} + \Y_w^{red} = K_{\Y/\Z} + \Y_w^{red}$$
\noindent and therefore that 
$$\O_\Y(D' + \Y_T - \Y_S^{red} + \Y_w^{red})|_{\Y_w^{red}} =
\omega_{\Y/\Z}(\Y_w^{red})|_{\Y_w^{red}}$$

\noindent The right hand side is equal to $\omega_{\Y_w^{red}}$ by the
adjunction formula.  To calculate the left hand side, we observe that
$w$ is not in $S$ by hypotheses, although it may be in $T$.  Thus there
exists an integer $m$ depending on the multiplicity of $w$ in $T$ for which
the following calculations hold:
\begin{eqnarray*} 
\O_\Y(D' + \Y_T - \Y_S^{red} + \Y_w^{red})|_{\Y_w^{red}} &=& \O_\Y(D' +
m\Y_w^{red}) |_{\Y_w^{red}}\\
&=&\O_\Y(D')|_{\Y_w^{red}} \otimes \O_\Y((\Y_w^{red})^{\otimes m})|_{\Y_w^{red}}\\
&\cong&\O_\Y(D')|_{\Y_w^{red}} \otimes \O_\Y \\
&\cong&\O_\Y(D')|_{\Y_w^{red}} \\
&=&\O_\Y(D' \cap \Y_w^{red})
\end{eqnarray*}

\noindent which proves Lemma \ref{lemma:canonical} and therefore Theorem
\ref{thm:irred}. $\blacksquare$

\begin{remark} \label{rmk:trivial} 

Note that we can identify $det(V^I)$ with a character of order 1 or 2
of ${\rm Pic}_{Weil}(\Y_v^{red})$.  Thus, $\epsilon(\Y_v^{red},
det(V^I)) = 1$, as it's the ratio of epsilon constants associated to
zeta functions.  In particular we can show that both $\epsilon(\Y_v,V)$
and $\epsilon(D',V)$ are trivial, which would give us another way of
proving Theorem \ref{thm:irred}.  However, for what follows it is more
illuminating to instead consider what their ratio is, and in particular
how close each is to being of the form $det(V^I)(K)$.

\end{remark}

\subsection{Partial Trivializations and the Canonical Cycles}

In this section we will describe in detail the relative canonical cycle
associated to line bundles with partial trivializations, as defined by
T. Saito in \cite{tS}, as well as other machinery which we will
need in order to compute the terms
$\epsilon(D'_v,V)^{-1}\epsilon(\Y_v,V)$ in the case where $\Y_v^{red}$
consists of more than one component.

\begin{definition} Let $\D$ be a divisor on a scheme $X$ and let
$\{\D_i\}_{i \in I}$ be the set of irreducible components of $\D$.  A
locally free sheaf $\E$ on $X$ is said to be partially trivialized on
$\D$ if there exists a family $\rho = (\rho_i)$ of $\O_{\D_i}$-morphisms
$\rho_i:\E|_{\D_i} \rightarrow \O_{\D_i}$ such that for all subsets $J
\subset I$, the map $\rho_J = \bigoplus_{i \in J} \rho_i: \E|_{\D_J}
\rightarrow \O_{\D_J}^J$  is surjective. \end{definition}

Given a partial trivialization of the sheaf $\E$ of rank $n$ on $X$,
Saito defines the relative top chern class $c_n(\E, \rho) \in H^{2n}(X$
$mod$ $\D, \Z_q(n))$ based on an idea of Anderson in \cite{A}. In
particular, Saito notes that there is a canonical isomorphism  $$\Phi:
H^{2n}(X mod \D, \Z_q(n)) \rightarrow H^{2n}(V mod \Delta,
\Z_q(n))$$ where $V$ is the covariant vector bundle associated to the
dual of $\E$.  We also have a natural map from $H^0(X, \Z_q)
\rightarrow H^{2n}(V$ $mod$ $\Delta, \Z_q(n))$. Let $[0]$ be the image of the
class $1 \in H^0(X,\Z_q)$ under this map, and then Saito defines  the
relative top chern class to be the inverse image of $[0]$ under the
canonical isomorphism $\Phi$ above.  Relative top chern classes satisfy nice
functorial properties, and the relative top chern class is mapped to
the normal top chern class under the canonical map $H^{2n}(X$ $mod$
$\D, \Z_q(n)) \rightarrow H^{2n}(X, \Z_q(n))$. Furthermore, the
following corollary of Proposition $1$ in \cite{tS} gives us a way to
compare the relative top chern classes associated to two different
partial trivializations.

\begin{corollary} \label{cor:compare} Let $X$ be a $\F_p$-scheme, and
let $(\E, \rho)$ be a partially trivialized locally free sheaf on $X$. 
Let $\sigma_i = f_i^{-1} \cdot \rho_i: \E |_{\D_i} \rightarrow \O_{\D_i}$
where $f_i$ comes from $\F_p^*$, so that $\sigma = (\sigma_i)$ is
another partial trivialization of $\E$.  Finally, let $\E_i =
Ker(\rho_i)$ so that $\rho|_{\D_i}$ is a partial trivialization of
$\E_i$.  Then we can compute the difference between the relative top
chern classes as $$c_n(\E, \rho) - c_n(\E, \sigma) = \sum \{f_i\} \cup
c_{n-1}(\E_i, \rho|_{\D_i})$$ \end{corollary}

In section 2 of \cite{tS}, Saito uses the construction of relative top
chern classes to define the relative canonical cycle.

\begin{definition} \label{defn:saito} Let $\D$ be a divisor with simple
normal crossings on a variety $X$ of dimension $n$ defined over a
perfect field $F$ of characteristic $p$, and let $U=X-\D$.  Let $\Omega^1_{X/F}(log$ $\D)$ be the locally free $\O_X$-module of rank $n$ of differential $1$-forms on $X$ with logarithmic poles along $\D$.  Then the cycle $$c_{X,U} = (-1)^n
c_n(\Omega^1_{X/F}(log\D), res)$$ is called the relative canonical cycle.  It lies inside the cohomology with compact support $H_c^{2n}(X mod$ $\D, \hat{\Z}'(n))$, where $\hat{\Z}' = \prod_{q \ne p} \Z_q$. The relative canonical cycle has degree equal
to $\chi(U_{\overline{F}}) = \Sigma (-1)^q dim H^q_c(U_{\overline{F}},
\overline{\Q}_l)$. Note that this definition differs from that of S.
Saito in \cite{sS}, but only up to a change in sign. \end{definition}

Saito observes that one can also define a relative top chern class (and
hence a relative canonical cycle) sitting inside of $H^n(X$ $mod$
$\D,\G_m)$, the divisor class group with modulus $D$,  and in
particular we can define $c_{X,U}$ as an element of $H^n(X$ $mod$
$\D,\G_m)$ in the case when $n=1$. For our work we will want to consider the case where $X$ is one of the components of $\Y_v^{red}$, and therefore is of dimension one. 
We will look at the relative top chern class $c_{X,U}$ lying
inside the generalized class group $$H^1(X mod \D, \G_m) = [(\oplus_{x
\notin \D}\Z) \oplus (\oplus_{x \in \D}K^*/U^1_x)]/K^*$$ \noindent
where $K$ is the fraction field of $X$ and $U^1_x = 1+m_x$. In
particular, the class $c_{X,U}$ can be computed in the following way
(see the example in \S 1 of \cite{tS}). Let  $\omega$ be
a nontrivial rational section of $\Omega^1_X(log$ $\D)$ such that for
all points $x \in \D$, $ord_x(\omega)= -1$ and $res_x(\omega) = 1$ then
the relative canonical cycle represented by the class of the zero cycle
which is supported off of $\D$ given by $$c_{X,U} = -\sum_{x \in
U}ord_x(\omega) \cdot [x]$$ 

\begin{proposition}
Let $\Y_v^{red}$ consist of two components $F'$ and $G'$. Let $D'$ be
a horizontal divisor chosen as in the previous sections.   

\begin{enumerate} 

\item[(1)]  There is a canonical isomorphism $\phi:  \O_{F'}(D' \cap F') \rightarrow \omega_{F'}(F' \cap G')$ up to multiplication by a global unit. 

\item[(2)] The global section $1 \in \Gamma(\O_{F'}(D' \cap F'))$ maps
under $\phi$ to an element $\gamma \in \Gamma(\omega_{F'}(F' \cap G'))$
such that $ord_x(\gamma) = 1$ if $x \in F' \cap D'$, $ord_x(\gamma) =
-1$ if $x \in F' \cap G'$, and  $ord_x(\gamma) = 0$ otherwise. 

\item[(3)] Set $a_x = res_x(\gamma)$ for all $x \in F' \cap G'$.
Then  $c_{F',U_{F'}}$ is such that $-c_{F',U_{F'}}$ is the class in $[(\oplus_{x \in (F' - G')}\Z) \oplus (\oplus_{x \in F' \cap G'}K^*/U^1_x)]/K^* = H^1(F'$ $mod$ $(F' \cap G'),
K)$ of the element

$$c= (\oplus_{x \in D' \cap F'} 1 \in \Z) \oplus (\oplus_{x \in F' - D'
- G'} 0 \in \Z) \oplus (\oplus_{x \in F' \cap G'}a_x)$$

\end{enumerate}
\end{proposition}

The proof of part $(1)$ follows from carrying through a series of
calculations analagous to those in the proof of Lemma
\ref{lemma:canonical}.  In particular,
\begin{eqnarray*}
\O_{F'}(D' \cap F') &=& \O_\Y(D') |_{F'} \\
&=&\O_\Y(D'+\Y_T) |_{F'}\\
&\cong&\omega_\Y(\Y_S^{red})|_{F'} \\
&=&\omega_\Y(\Y_v^{red})|_{F'} \\
&=& [\omega_\Y(F') \otimes \O_\Y(G')] |_{F'}\\
&=& \omega_\Y(F')|_{F'} \otimes \O_\Y(G') |_{F'}\\
&=& \omega_{F'} \otimes \O_{F'}(F' \cap G') \\
&=& \omega_{F'}(F' \cap G')
\end{eqnarray*}

To prove parts $(2)$ and $(3)$ of the proposition, we set $X = F'$ and $\D = F' \cap G'$.  Proving these statements is then just a matter of calculating the
various orders and residues of $\gamma$ given that we know them for the
element  $1 \in \Gamma(\O_{F'}(D' \cap F'))$. Explicitly, they can be
computed by following the residue map on elements of the sheaves
through the equalities and congruences in the calculations above.  Note
that all of the isomorphisms are unique with the exception of
$\O_\Y(D'+\Y_T) \cong \omega_\Y(\Y_S^{red})$.  This map, while not unique, is well-defined up to multiplication by a global unit, and
therefore when we look at classes mod $K^*$  the discrepancy will not
matter.

This proposition gives us an explicit way to construct the relative
canonical class in our situation. In particular, the $a_x$ terms come
about because of the difference in natural partial trivializations on the
sheaves $\O_{F'}(D' \cap F')$ and $\omega_{F'}(F' \cap G')$ associated
to the restriction map $\O_{F'}(D' \cap F') \rightarrow \O_{F'}(D' \cap
F')|_{F' \cap G'} = \O_{F' \cap G'}$ and the residue maps $res_x$.

For the computations in the next section, we will need the following
definitions.

\begin{definition}  \label{defn:classes}  Define the following classes
which lie in the generalized class group $H^1(F'$ $mod$ $(F' \cap G'),
K)$.

\begin{enumerate}

\item[a.] Define an element $\lambda \in (\oplus_{x \in (F' - G')}\Z)
\oplus (\oplus_{x \in F' \cap G'}K^*/U^1_x)$ which has components equal
to $1 \in \Z$ at all points $x \in D' \cap F'$, equal to $0 \in \Z$ at
all points $x$ in $F' - D' - G'$ and equal to the identity in
$K^*/U^1_x$ for all points $x \in F' \cap G'$. We then look at the
class $[\lambda]_{F'} \in H^1(F'$ $mod$ $(F' \cap G'), K)$, which is
the first relative chern class of the line bundle $\O_{F'}(D' \cap F')$
with partial trivializations.  One can define $[\lambda]_{G'}$ in a similar way.

\item[b.] Let $\delta_{F'}$ be the class in $H^1(F'$ $mod$ $(F' \cap
G'), K)$ which corresponds to the element $\delta=(\oplus 0) \oplus
(\oplus a_x) \in (\oplus_{x \in F' - G'}\Z) \oplus (\oplus_{x F' \cap
G'}K^*/U^1_x)$.  In other words, this element is trivial at all places
corresponding to $x \notin F' \cap G'$ and for those places which
correspond to points $x \in F' \cap G'$ consists of the  terms $a_x$
coming about as the difference between the partial trivializations of
$\O_\Y(D' + \Y_T)$ and $\omega_\Y(\Y_S^{red})$, as found in the above
characterization of $c_{F',U_{F'}}$.  Note that $\delta_{F'}$ can be thought
of as the quotient of $c_{F',U_{F'}}$ and $[\lambda]_{F'}$. 
One can define $\delta_{G'}$ in a similar way.

\end{enumerate}
\end{definition}

\subsection{The General Case}

We have shown that in the situation where $\Y_v^{red}$ consists of a
single component $F'$ then the fibral contribution to the root number
is positive. Now we will consider the next case, where $\Y_v^{red}$
consists of two irreducible components, say $F'$ and $G'$.  Note that
in particular this implies that $v \in S$. Recall that from Equation
$3.1$ above we are interested in comparing $\epsilon(D_v',V)$ and
$\epsilon(\Y_v,V)$.  By our initial assumptions, $D'$ intersects
$\Y_v^{red}$ in smooth points of $\Y_v^{red}$, so in particular we get
that the set $D' \cap F' \cap G' = \emptyset$

Define $I_1=I_{\mu(F)}$ and $I_2=I_{\mu(G)}$, where $F$ and $G$ are
components of $\X_v^{red}$ lying above $F'$ and $G'$ respectively. Then
$det(V^{I_1})$ is a character of the Galois group of the cover $F
\rightarrow F'$, which will be tame with respect to the divisor $F'
\cap G'$.  Classfield theory says that we can therefore view
$det(V^{I_1})$ as a character of the ray class group of $F'$ with
conductor $F' \cap G'$. We wish to define the term $det(V^{I_1})
(\pi_{D',y})$ for points $y \in F'$.  In order to do so, we view $det(V^{I_1})$ as a
character of the ideles $J_{F'}$ of $F'$.  In other words, it is an
idele class character modulo the conductor, which will be supported on
$F' \cap G'$.  We then define $det(V^{I_1})(\pi_{D',y})$ to be the value
of $det(V^{I_1})$ on the idele $(1,\ldots,1,\pi_{D',y},1,\ldots)$ which is
trivial away from $y$.  This is well defined as the conductor of
$det(V^{I_1})$ does not involve $y$ and the difference between two
local uniformizers is a unit.

If we define $det(V^{I_1})(D' \cap F')$ to be equal to the product
$\prod_{y \in D' \cap F'} det(V^{I_1})(\pi_{D',y})$, then this term will
be independent of the choices of uniformizers as all components are
unramified, and we are able to make the following calculation:
\begin{eqnarray}
\epsilon(D'_v, V) &=& \prod_{y \in D' \cap \Y_v^{red}}
\epsilon(y,V)\nonumber\\
&=&\prod_{y \in D' \cap F'} \epsilon(y,V)\prod_{y \in D' \cap G'}
\epsilon(y,V)\nonumber\\
&=&\prod_{y \in D' \cap F'} det(V^{I_1})(\pi_{D',y})\prod_{y \in D'
\cap G'}det(V^{I_2})(\pi_{D',y})\nonumber\\
&=&det(V^{I_1})(D' \cap F') \cdot det(V^{I_2})(D' \cap G')
\end{eqnarray}

Recall that in the case where $\Y_v^{red}$ consisted of a single
component $F'$, we were able to show that $\epsilon(D'_v,V) =
det(V^I)(D' \cap F')$. In that case Lemma \ref{lemma:canonical} showed
that our hypothesis on $D'$ implied that $D' \cap F'$ was a canonical
divisor on $F'$. The preceding section showed that in this more
complicated case, while $D' \cap F'$ is not a canonical divisor on
$F'$, it is close to being one.  To make this precise requires the
results of the previous section. In particular, when viewed as an idele
class character, $det(V^{I_1})$ breaks into components $det(V^{I_1})_x$
which are unramified for all $x \notin F' \cap G'$, and therefore we get
$$det(V^{I_1})(D' \cap F') = \prod_{y \in D' \cap F'}
det(V^{I_1})_y(\pi_{D',y}) = det(V^{I_1})([\lambda]_{F'})$$

Therefore Equation $3.2$ says that in the case where $V$ is an orthogonal virtual representation of
dimension $0$ and trivial determinant, $\Y_v^{red}$ consists of two
components $F'$ and $G'$ and $D'$ is chosen as above, then we have
that 
\begin{eqnarray}
\epsilon(D'_v,V) &=& det(V^{I_1})([\lambda]_{F'})
det(V^{I_2})([\lambda]_{G'}) 
\end{eqnarray}
Switching gears, we now want to take a look at the term
$\epsilon(\Y_v,V)$. For the moment, we will assume that $F' \cap G'$
consists of a single point $z$. We begin by looking at the two exact
sequences:
$$ 0 \rightarrow U = \Y_v^{red} - z \rightarrow \Y_v^{red} \rightarrow z
\rightarrow 0 $$
$$ 0 \rightarrow U \rightarrow \widehat{\Y_v^{red}} = F' \amalg G'
\rightarrow \{z_{F'}, z_{G'} \} \rightarrow 0 $$

\noindent where $z_{F'}$ (respectively $z_{G'}$) is the point $z$ thought of
as sitting just on $F'$ (respectively $G'$). Epsilon factors are multiplicative within exact sequences as well as in disjoint unions, so these sequences imply that
\begin{eqnarray}
\epsilon(\Y_v^{red},V) &=& \epsilon(U,V) \epsilon(z,V) \nonumber\\
&=& \frac{\epsilon(F',V)
\epsilon(G',V)}{\epsilon(z_{F'},V)\epsilon(z_{G'},V)} \epsilon(z,V)
\end{eqnarray}

To continue, we must consider the $\epsilon(F',V)$ term. In order to
compute this term we use the following result proven by Saito in
\cite{tS}

\begin{lemma} 
\label{lemma:saito} 

Let $X,U$ be as in Definition \ref{defn:saito}, and let the action of
$G$ be \'etale on $U$.  Then $\prod_{y \in U}\epsilon_y(X,V) =
det(V)(c_{X,U})$ 

\end{lemma}

Applying this lemma to our situation, we are able to make the following
computation:
\begin{eqnarray}
\epsilon(F',V) &=& \epsilon(F',V^{I_1})\nonumber\\
&=& \prod_{y \in (F')^0} \epsilon_y(F',V^{I_1}) \nonumber\\
&=& \epsilon_z(F',V^{I_1}) \prod_{y \neq z} (\epsilon_y(F',V^{I_1}))
\nonumber\\
&=& \epsilon_z(F',V^{I_1})det(V^{I_1})(c_{F',U_{F'}}) \nonumber\\
&=&\epsilon_{0,z}(F',V^{I_1}) \epsilon(z_{F'},V^{I_1}) det(V^{I_1})(c_{F',U_{F'}}) 
\end{eqnarray}

Plugging Equation $3.5$ (and the analogous formula for $\epsilon(G',V)$)
into Equation $3.4$ gives that 
\begin{eqnarray*}
\epsilon(\Y_v^{red},V)&=&\epsilon_{0,z}(F',V^{I_1})
\epsilon_{0,z}(G',V^{I_2})
det(V^{I_1})(c_{F',U_{F'}})\det(V^{I_2})(c_{G',U_{G'}})\epsilon(z,V)
\end{eqnarray*}

\noindent which we we can combine with Equation $3.3$ to get that
$$\frac{\epsilon(\Y_v,V)}{\epsilon(D'_v,V)}
=\frac{det(V^{I_1})(c_{F',U_{F'}})\det(V^{I_2})(c_{G',U_{G'}})}
{det(V^{I_1})([\lambda]_{F'}) det(V^{I_2})([\lambda]_{G'}) }
\epsilon_{0,z}(F',V^{I_1})\epsilon_{0,z}(G',V^{I_2})\epsilon(z,V)$$

\noindent Note that
$$\frac{det(V^{I_1})(c_{F',U_{F'}})}{det(V^{I_1})([\lambda]_{F'})} =
det(V^{I_1})(\delta_{F'})$$

\noindent where $\delta$ is the class defined in Definition
\ref{defn:classes}.

Considering a slightly more general case, in which we still only have
two components, but where $F' \cap G'$ consists of more than one point,
it is clear that all of the calculations will follow through and we
will get that 
$$\frac{\epsilon(Y_v,V)}{\epsilon(D'_v,V)} =
det(V^{I_1})(\delta_{F'})det(V^{I_2})(\delta_{G'})\prod_{z \in F' \cap G'}
\epsilon_{0,z} (F',V^{I_1}) \epsilon_{0,z} (G',V^{I_1}) \epsilon(z,V)$$

If we have more than two components in $\Y_v^{red}$ then the
bookkeeping becomes more complicated but the mathematics does not.  We
first set up the necessary notation.  Let $C_i$ be the components of
$\Y_v^{red}$. Furthermore, let $C_{i,j} = C_i \cap C_j$, $Z = \cup_{i
\ne j} C_{i,j}$ be the collection of all intersection points and let
$U_{C_i}$ be the open set consisting of $C_i - Z$. Finally, let $I_i$
be the inertia group associated to $C_i$ as above.  We are still
interested in computing $\epsilon(\Y_v,V)$ and $\epsilon(D'_v,V)$.  Let
$\lambda_v$ and $\delta_{v,C_i}$ be the classes $\lambda$ and $\delta$
defined above for a particular class $v$ and a particluar component
$C_i$.  In particular, recall that $\delta_{v,C_i}$ can be calculated purely
from looking at points $z \in Z$

For the latter, the computation works just as it did before, as we know
that if $i < j$ the $C_{i,j}$ are disjoint from each other as well as
from $D'$.  We obtain that

\begin{eqnarray*}
\epsilon(D'_v,V) &=&\prod_i det(V^{I_i})(C_i \cap D') \\
&=& \prod_i det(V^{I_i})([\lambda_v]_{C_i})
\end{eqnarray*}

To compute $\epsilon(\Y_v,V)$ we need to use the following exact
sequences:
$$ 0 \rightarrow U = \Y_v^{red} - Z \rightarrow \Y_v^{red} \rightarrow
Z \rightarrow 0 $$
$$ 0 \rightarrow U \rightarrow \widehat{\Y_v^{red}} = \amalg_i C_i
\rightarrow \amalg_{i \ne j} C_{i,j} \rightarrow 0 $$
\noindent where $\amalg_{i \ne j} C_{i,j}$ can be thought of as the set
consisting of two copies of $Z$, with each point considered as sitting
once on each of the two $C_i$ which it comes from originally. We can
now use these sequences as well as the above calculations of
$\epsilon(C_i,V)$ to get that
\begin{eqnarray*}
\epsilon(\Y_v^{red},V) &=& \epsilon(U,V) \epsilon(Z,V) \\
&=& \prod_i \epsilon(C_i,V) \prod_{z \in Z} \frac{\epsilon(z,V)}
{\epsilon(z_{C_{i_1}},V)\epsilon(z_{C_{i_2}},V)} \\
&=& \prod_i det(V^{I_i})(c_{C_i,U_{C_i}}) \prod_{z \in Z}
\epsilon_{0,z}(C_{i_1},V^{I_{i_1}})\epsilon_{0,z}(C_{i_2},V^{I_{i_2}})
\end{eqnarray*}
\noindent where we think of $z \in Z$ as lying on $C_{i_1} \cap
C_{i_2}$. If we put all of these calculations together we get the
following result.

\begin{theorem} 
\label{thm:prodform} 

Under all of the above hypotheses and notation, we get that for all
$v$, $$\frac{\epsilon(\Y_v,V)}{\epsilon(D'_v,V)}
=\prod_i det(V^{I_i})(\delta_{v,C_i})\prod_{z \in Z}
\epsilon_{0,z}(C_{i_1},V^{I_{i_1}})\epsilon_{0,z}(C_{i_2},V^{I_{i_2}})
\epsilon(z,V)$$

\noindent where both of these products are equal to one if the set $Z$
is empty.

\end{theorem}

Combining Equation (\ref{eq:doit}), Lemma \ref{lemma:fq} and Theorem
\ref{thm:prodform} gives a precise form of Theorem \ref{thm:main}. 
Note that other than the term $\epsilon_\infty(\Y,V)$, the other terms
depend only on the crossing points of the components of fibers
$\Y_v^{red}$.

\section{Examples}
\label{section:examples}
\setcounter{equation}{0}

In this section we wish to show several examples that apply Theorem
\ref{thm:main}.  In order to do this, we must first have concrete examples of
finite groups acting tamely on arithmetic surfaces. We find one class of such
examples by using the following result from the thesis of Seon-In Kwon
{\cite{K}}

\begin{theorem}
\label{thm:elliptic}

Let $X$ be an elliptic curve over $K$ and let $\X$ be the minimal model
of $X$ over $\O_K$. Consider the action of a group $G\cong \Z/n\Z
\times \Z/m\Z$ $\subset X(K)$ of torsion points on $\X$. Then the
action of $G$ on $\X$ is numerically tame if and only if for each place
$v$ of $\O_K$ whose residue characteristic $p$ divides the order of
$G$, the following conditions are satisfied:

\noindent (i) The minimal model $\X$ has good or multiplicative
reduction at $v$.

\noindent (ii) The Zariski closure in $\X$ of the $p$-Sylow
subgroup $G_p$ of $G$ is smooth over $Spec \O_K$.

In particular, these conditions imply that $gcd(n, m)=1$.
\end{theorem}

This theorem provides us with a set of concrete criteria for checking when a
finite group acts tamely on the integral model of an elliptic curve as
numerical tameness implies tameness.  In particular, the second condition asks
us to compute the $p$-torsion points of the minimal model over $p$, and check
that they do not coalesce when we reduce mod $p$. 

We will now show an example of a computation of the orthogonal
$\epsilon$-constants associated to the tame action of a finite group on a
surface. First, we must calculate terms $\epsilon(z,V)$, where $z \in \Y$ is a
closed point defined over a finite field.  Section $2.5$ of \cite{CEPT1} gives
us the following way of making this computation.

\begin{lemma} 
\label{lemma:points} 

Let $x$ be a point of $\X$ over a point $y \in \Y$ which has finite residue
field.  Furthermore, let $F_{x}$ be the arithmetic Frobenius element lying in
$G$.  Then $\epsilon(y,V) = det(V^{I_x})(-F_{x})$, where $I_x \subseteq G$ is
the inertia group of the point $x$.

\end{lemma}

Let $X$ be the elliptic curve given by the equation $y^2+xy+y = x^3$.  This
equation is minimal over every prime $p \in Spec(\Z)$, and thus it defines
$\X$, the minimal model over $\Z$.  The torsion subgroup of $X$ is isomorphic
to $\Z/3\Z$, and the torsion points of order three are $(0,0)$ and $(0,-1)$. 
We wish to check to see whether or not the action of $G \cong \Z/3\Z$ is
numerically tame by the criteria in Theorem \ref{thm:elliptic}.  In order to do
this, we first note that the discriminant of $\X$ is $-26 = -1 \cdot 2 \cdot
13$, and thus $\X$ has good reduction at $3$.  Furthermore $\X$ has
multiplicative reduction at $2$ and at $13$ (and in particular the fibers have
Kodaira type $I1$).

Next we check condition (ii).  $G_3 = G = \{(0,0),(0,-1),{\bf 0}\}$, and these
points clearly do not coalesce mod $q$ for any prime $q$ (and in particular for
$q=3$).  Thus, this action of $G$ on $\X$ does satisfy the appropriate
conditions.

Therefore the action of $G$ is in fact numerically tame on $\X$.  Furthermore,
it follows from formulae of Velu in {\cite{V}} that $\Y = \X /G$ is the
integral model of the elliptic curve defined by the equation $y^2+xy+y =
x^3-5x-8$.   The fibers of $\Y$ are also nonsingular with the exceptions of the
fibers at $v=2, 12$ which are both of Kodaira type $I3$.  However, $\Y$ is not
regular, and thus the results of Section \ref{section:main} do not apply and in
fact the $\epsilon$-constants are not well-defined.  However, due to Theorem
$3.8$ of Kwon's dissertation \cite{K}, we know that after a finite number of
blow-ups on the singular fibers we can blow up $\X$ in a way such that the
action of the group $G$ extends to a tame action of $G$ on the blow-up $\X_1$,
and the quotient $\Y_1 = \X_1/G$ is in fact regular.  This theorem applies
because all of the local decomposition groups must be subgroups of $\Z/3\Z$ and
in particular must be cyclic of degree $n \le 3$.

Next we need to define a representation $V$ of $\Z/3\Z$ satisfying certain
properties.  We know from representation theory that there are two distinct
nontrivial one-dimensional characters of $\Z/3\Z$ of order three.  Let us
define $V_1$ to be the sum of these characters and $V_2$ to be $2\chi_0$, where
$\chi_0$ is the trivial character.  We then define $V$ to be $V_1 - V_2$.  Sums
of characters are always orthogonal representations, so $V$ will be
orthogonal.  It also is not hard to see that $V$ has dimension zero and trivial
determinant.

In general, computing $\epsilon(\Y_1,V)$ might be difficult, but in light of
Theorem \ref{thm:main}, the computation simplifies greatly. In particular, we
only need to compute $\epsilon_{\infty,0}(\Y_1,V)$,
$det(V^{I_j})(\delta_{v,C_i})$ for $v = 2, 13$, and the terms
$$\epsilon_{0,z}(C_{i_1},V^{I_{i_1}})\epsilon_{0,z}(C_{i_2},V^{I_{i_2}})
\epsilon(z,V)$$ at the singular points above the primes $p=2$ and $p=13$.  For
the above choice of the representation $V$, we can see that $det(V^I)$ is
trivial for all possible inertia groups $I$.  More precisely, $det(V_j^I)$ will
be trivial for $j=1,2$.  If $I$ acts trivially on $V_j$ this is obvious, as the
$det(V_j)$ are both in fact trivial.  On the other hand, if $I$ acts
nontrivially on $V_j$, then $V_j^I$ will be trivial as the kernels of both
characters which make up $V_j$ are the same, and thus $det(V_j^I)$ will be
trivial as well.

Let us first look at the part of the calculation of $\epsilon(\Y_1,V)$ coming
from the fiber of $\Y_1$ above the prime $2$.  Denote the three components of
$\Y_2$ by $F_1, F_2, and F_3$.  Let $I_i$ be the inertia group associated to
$F_i$.  In particular, $det(V^{I_i})$ is trivial in each of these cases for the
reasons described above.  Thus, the $det(V^{I_j})(\delta_{2,C_i})$ terms are
equal to one.  Many of the $\epsilon_{0,z}(C_{i},V^{I_{i}})$ terms will also
immediately be equal to one as many of the $V^{I_i}$ terms are themselves
trivial.  To compute the others, we use Theorem $2$ of Saito in \cite{tS}. 
Because we are looking at cases where $det(V^I_i)$ is trivial, these formulae
reduce the computation of $\epsilon_{0,z}(C_{i},V^{I_{i}})$ to the computation
of a Gauss sum $\tau_{C_i}(V^{I_i})$.  We can now use the fact that our
representation $V$ is the sum of a representation and its complex conjugate to
get that  $\tau_{C_i}(V^{I_i}) = 1$.

The above paragraph holds for the points above $p=13$ as well, so we can ignore those
terms. We can now use Lemma \ref{lemma:points} in order to compute the
$\epsilon(z,V)$ terms.  In particular, the fact that all of the $det(V^I)$ terms are trivial tells us that these terms are also equal to one.  To summarize, we have that $\epsilon(\Y_1,V) = \epsilon_{\infty,0}(\Y_1,V)$.  

In Theorem 4.0.1 of \cite{CPT} they show that the Euler characteristics and the character functions $\zeta_S$ associated to the action of a finite group on a minimal model over $\Z$ of an elliptic curve over $\Q$ satisfying certain properties (which our example does satisfy) are trivial.  This is because the group must act trivially on the various $H^{p,q}$ pieces of the Hodge structure.  But this in turn shows that $\epsilon_{\infty,0}(\Y_1,V)$ is trivial, and thus $\epsilon(\Y_1,V) = 1$.

We note that many of the computations in the last example will hold whenever we
are in the case of one of Kwon's examples.  In particular, it will often be the
case that we can show that $det(V^{I_i})$ is trivial.  Combining this with the
above-cited results in \cite{CPT} shows that $\epsilon(\Y,V)$ is trivial for a
large number of examples where $\Y$ is a blowup of the minimal model of an
elliptic curve with $G$ acting tamely as in \cite{K}.

\end{document}